\documentclass[12pt]{msml2020} 

\title[Quasi-optimal DNN]{Analysis of Deep Neural Networks with Quasi-optimal polynomial approximation rates}
\usepackage{times}
\usepackage{amssymb}
\usepackage{bm}

\newtheorem{assumption}{Assumption}

\msmlauthor{\Name{Joseph Daws} \Email{jdaws@vols.utk.edu}\and
\Name{Clayton Webster} \Email{cwebst13@utk.edu}\\
\addr Department of Mathematics, The University of Tennessee, Knoxville, TN 37996, USA}

\begin{document}

\maketitle

\begin{abstract}%
	We show the existence of a deep neural network capable of approximating a wide
	class of high-dimensional approximations. The construction of the proposed neural network 
	is based on quasi-optimal polynomial approximation.
	We show that this network achieves an error rate that is sub-exponential 
	in the number of polynomial functions, $M$, used in the polynomial approximation. 
	The complexity of the network which achieves this sub-exponential rate is shown to be algebraic in $M$.
\end{abstract}

\begin{keywords}%
  Expressivity, Deep Neural Networks, Approximation Theory
\end{keywords}

\section{Introduction}
We propose a class of neural networks whose architecture is inspired by high-dimensional 
tensor products of polynomials. Our main result shows there exists a deep neural network (DNN) 
which achieves the same error rate as a quasi-optimal $M$-term polynomial approximations 
for approximating a given target function $u$ satisfying some mild assumptions which are described below.
The proof of our main theorem is constructive in the sense that we define
the values of all parameters of the approximating network.  
The architecture of our proposed class of networks is similar to those presented in 
\cite{deep_learning_schwab,new_schwab,Yarotsky_2017,deep_relu_du} 
whose architectures are chosen based on polynomials when the coefficient values are exactly known.
The analysis in these works uses exact coefficient values, which may be difficult to obtain
for many problems of interest. In contrast, the analysis presented here relies on a 
bound for the coefficients which satisfies Assumption \ref{eq:bounds_assumption}.
The improved error rates and reduce complexity presented in this work depend on
the constructive polynomial \textit{quasi-optimal approximation}. In paricular,
the error rate is very similar to the one obtained in \cite{analysis_of_quasi_optimal}. 
The analysis of quasi-optimal approximations presented in that work applies to
functions which satisfy some assumptions on its regularity.
\begin{assumption} \label{as:function_of_interest}
	The function $u: \mathbb{R}^d \rightarrow \mathbb{R}$ satisfies 
	\begin{enumerate}
		\item $u$ is analytic in a poly-ellipse.
		\item when expanded in a polynomial basis or frame $\{ \Psi_{\pmb{\nu}} \}$,
		         where each $\Psi_{\pmb{\nu}}: \mathbb{R}^d \rightarrow \mathbb{R}$ is a tensor product of 
		         one-dimensional polynomials, the coefficients $c_{\pmb{\nu}}$ are bounded 
		         exponentially
		         	\begin{equation}
					\| c_{\pmb{\nu}} \| \le e^{-b(\pmb{\nu})}
				\end{equation}
			for some function $b(\pmb{\nu})$ which satisfies \cite[Assumption~3]{analysis_of_quasi_optimal}
			which is stated as Assumption \ref{eq:bounds_assumption} in section \ref{sec:quasi_optimal} of this work.
	\end{enumerate}
\end{assumption}
The work \cite{analysis_of_quasi_optimal} shows that the solutions to a wide class of elliptic PDE 
satisfy Assumption \ref{as:function_of_interest}. However, the class of functions which satisfy this
assumption is much broader than solutions to certain PDE and is a mild assumption to make
for functions with some expected smoothness. 

Suppose $u$ is a function of interest which satisfies Assumption \ref{as:function_of_interest}. 
Let $u_{\Lambda_M}$ be the polynomial approximation given by 
\begin{equation} \label{eq:lambda_truncation}
	u_{\Lambda_M}(\bm{y}) = \sum_{\bm{\nu} \in \Lambda_M} c_{\bm{\nu}} \Psi_{\bm{\nu}} (\bm{y}),
\end{equation}
where $\bm{y} \in \mathbb{C}^d$, $\Lambda_M \subset \{ \bm{\nu} = (\nu_i)_{1\le i \le d} : \nu \in \mathbb{N}\}$
such that $\#(\Lambda_M) \le M$ and $c_{\bm{\nu}}$ 
are the coefficients associated with the polynomial $\Psi_{\bm{\nu}}$ which we will assume is a tensor 
product of polynomials.
The error between $u$ and $u_{\Lambda_M}$ measured by the quantity
	\begin{equation} \label{eq:best_M_term_error}
		\| u - u_{\Lambda_M} \|,
	\end{equation}
for an appropriate choice of norm. This quantity is of keen interest
in approximation theory. In particular, when $\Lambda_M$ is 
the set of indices associated with the $M$ largest values of $c_{\Lambda_M}$,
\eqref{eq:best_M_term_error} is called the \textit{best $M$-term error}. Bounds on 
the best $M$-term error have been extensively studied for many different bases and contexts,
see e.g. \cite{rons_good_book}. The construction of an $M$-term 
approximation is trivial once the index set $\Lambda_M^{best}$ associated to the
best $M$-term error is identified. In general, a brute force computation of of such 
index sets in untenable. Algorithms for the construction of 
close approximations of $u_{\Lambda_M^{best}}$ have been proposed. 
For example, see \cite[Section~8]{rons_good_nonlinear_approximation_review} 
which reviews several such algorithms.

In this work, we chose $\Lambda_M$ to be the quasi-optimal index set 
as derived in the work \cite{analysis_of_quasi_optimal}. 
In the quasi-optimal framework, the indices in $\Lambda_M$ are chosen based on sharp bounds
of the coefficients $c_{\bm{\nu}}$ in \eqref{eq:lambda_truncation} rather than exact values of 
coefficients which are often more expensive to compute. The architecture of our network will be 
explicitly parameterized by the quasi-optimal index set. 
Moreover, we will assume that it has already been identified so
that we may establish a quasi best $M$-term network error estimate for approximation by a neural network. 
Such an error estimate is useful for bench-marking the performance of other neural networks which have 
comparable complexity to the one we construct. 

In order to quantify the error and the complexity of a network it is convenient to identify it both with a function, 
i.e., a mapping from the $d$-dimensional input space to a real value, and a graph, i.e., an acyclic graph 
which describes the arrangement and connections of neurons in the network.
For a full description of these interpretations see Section \ref{sec:nn}. 
Based on context it should be obvious which interpretation is used.
It is clear from its functional interpretation that the error of the network can be q
uantified by considering
	\begin{equation} \label{eq:network_error}
		\| u - u_{NN} \|
	\end{equation}
for any desired norm where $u_{NN}$ is the approximating network which we will
construct. 

The quantity \eqref{eq:network_error} is difficult to analyze without intimate 
knowledge of the parameters associated with $u_{NN}$. 
Following a constructive approach similar to the work 
\cite{deep_learning_schwab}, we will leverage the approximation power of polynomials
to both inform the architecture of the network $u_{NN}$ and bound the quantity \eqref{eq:network_error}.
By the triangle inequality, for any norm it is clear that
	\begin{equation} \label{eq:tricky_networks}
		\| u - u_{NN} \|  \le \| u - u_{\Lambda_M} \| + \| u_{\Lambda_M} - u_{NN} \|.
	\end{equation}
The approximation power of quasi-optimal approximations is established \cite{analysis_of_quasi_optimal}.
Therefore, the first term can be bounded by a sub-exponential expression in $M$. 
The main task is to design a neural network which can approximate 
$u_{\Lambda_M}$ with arbitrary accuracy. Finally, by choosing the
error between $u_{\Lambda_M}$ and $u_{NN}$ to be on the same order as the first term in 
\eqref{eq:tricky_networks} we will show that our network achieves the same sub-exponential 
rate as the quasi-optimal approximate. We will provide an explicit construction 
of a network $u_{NN}$ that approximates the quasi-optimal polynomial $u_{\Lambda_M}$
arbitrarily well in Section \ref{sec:nn}. 

\subsection{Context and Related Works}
According to \cite{Poggio_2017}, theoretical concerns about neural networks fall in three broad categories:
\begin{enumerate}
	\item Expressibility -- Given a network architecture what kinds of function mappings is it capable of approximating?
	\item Optimization -- How to identify parameters associated with a network so that it achieves the desired task?
	\item Generalization -- Once one a set of parameters associated to a network are fixed by a training process, 
	                                     how well does a network perform the task on data not in the training set, i.e., the testing data?
\end{enumerate}
These questions, at least superficially, have a vague resemblance to the central quesitons
that approximation theory attempts to answer. In order to contextualize the utility of our work, 
we will now draw analogues between mathematical approximation theory and neural networks
as an approximation paradigm. Many classical approximation theory results 
consider first a large space of where an object of interest lies. 
Often our object of interest is unattainable and therefore a subspace, called the approximation space, 
is chosen. Theorems are then proved showing how close the best approximation 
from this approximation space is to the desired function of interest. Such theorems are often termed 
\textit{direct theorems}. Our main result, Theorem \ref{thm:main_thm}, 
is a direct theorem about approximation by a class of neural networks. 
Continuing our analogue to approximation theory, the issue of identifying parameters in a network is clearly
analogous to constructing an algorithm which can construct or at least approximate the best approximation given
theoretically from a direct theorem. 
In short, the object of interest, the approximation space and the algorithm used to 
construct an approximation are fundamentally related.
Therefore, although direct theorems do not explicit provide or construct an approximation, 
they are critical for evaluating how well a constructive algorithm performs.

The power of neural networks is evident by their successful deployment in solving
extremely challenging problems such as image classification \cite{MR2224485, NIPS2012_4824},
artificial intelligence in playing games \cite{alpha_go}, novel data generation of the same type as the 
inputs \cite{NIPS2014_5423} and prediction of time series 
evolution \cite{LSTM}. Our work is one of several that show that neural networks
have immense expressive power. Shallow networks were shown to be able to approximate
any continuous function by Cybenko \cite{Cybenko1989ApproximationBS}. However,
recent works has shown that using deep networks allows one to achieve the same 
expressive power with fewer overall free parameters \cite{telgarksy_benefit_depth,Poggio_2017}.

The expressive power of deep neural networks (DNN) in the context of function approximation 
has been examined previously \cite{NIPS2017_7203}.
These works however offer no suggestion of the best way to add complexity to a given network in order to increase
its accuracy to the desired function. Since our neural network is based on quasi-optimal polynomial approximations
it is clear how one could increase accuracy. We need to simply add more indices to the quasi-optimal set and therefore
add more polynomial blocks to our network. Such an addition can be performed to a network which already fully trained
without any need to adjust the trained parameters any further. 

\subsection{Outline of Paper}
The constructive proof of our proposed networks relies on
many details related to polynomial and in particular quasi-optimal approximation.
We review the relevant theory in Section \ref{sec:quasi_optimal}. In Section 
\ref{sec:nn} we review some basic concepts about neural networks. Then we 
recall several recent works for approximating polynomials with neural networks.
Our main approximation theorem is proven in Section \ref{sec:quasi_net}.

\section{Quasi-optimal polynomial approximations} \label{sec:quasi_optimal}

In this section, we briefly review the results of \cite{analysis_of_quasi_optimal} which 
analyzes quasi-optimal index sets associated with high-dimensional tensor product polynomial 
approximations and we recall the sub-exponential convergence rate,
	\begin{equation}
		M \exp(-(\kappa M)^{1/d})
	\end{equation}
where $M$ is the number of terms used in the expansion and 
$\kappa$ is a constant which doesn't depend on the dimension $d$.
One may find sharp bounds for $c_{\bm{\nu}}$ by either \textit{a priori} or \textit{a posteriori} 
means. Once the bounds are established the quasi-optimal index $\Lambda_M$ is chosen corresponding to the
$M$ largest bounds. It has also been shown that the quasi-optimal method performs similarly to best $M$-term methods
\cite{quasi-optimal-stochastic-galerkin, optimal-stochastic-galerkin} while having a reduced cost since
exact coefficients are not computed.

\subsection{Coefficient bounds}
The error between any $M$-term expansion in a given set of functions 
$\{ \Psi_{\pmb{\nu}} \}$ associated with an index set $\Lambda_M$
and the function $u$ involves a sum of the norms of the coefficients with indices not in $\Lambda_M$.
The error of this approximation is given by the quantity $\| u - u_{\Lambda_M} \|$ with 
some appropriately chosen norm. For any set $\Lambda_M$ of $M$ multi-indices,
	\begin{equation}
		\| u - u_{\Lambda_M} \| =
		\left \| \sum_{\pmb{\nu} \in \Lambda_M \cup \Lambda_M^c} 
		c_{\pmb{\nu}} \Psi_{\pmb{\nu}} 
		- \sum_{\pmb{\nu} \in \Lambda_M} c_{\pmb{\nu}} \Psi_{\pmb{\nu}} \right \|
		\le
		\sum_{\pmb{\nu} \in \Lambda_M^c} \left \|  c_{\pmb{\nu}} \right \|.
	\end{equation}
Hence, the approximation $u_{\Lambda_M}$ is optimal if it is chosen to contain
the indices associated to the largest $M$ values of $\| c_{\pmb{\nu}} \|$. 
As previously mentioned, it is often difficult to construct such an approximation.
Alternatively, a bound for each of the coefficients can also be used to choose an index set.
Consider the case when each coefficient is bounded by a function $B(\pmb{\nu})$ then  
for any $\Lambda_M$,
	\begin{equation} \label{eq:upsilon_error}
		\| u - u_{\Lambda_M} \| \le \sum_{\pmb{\nu} \in \Lambda_M^c} \left \|  c_{\pmb{\nu}} \right \|
		\le \sum_{\pmb{\nu} \in \Lambda_M^c} B(\pmb{\nu}).
	\end{equation}
Therefore, an approximation can be constructed by choosing the indices corresponding to the 
$M$ largest values of $B(\pmb{\nu})$. Hence, if 
$B(\pmb{\nu})$ is known for a set of functions $\{ \Psi_{\pmb{\nu}} \}$ and is a
reasonably sharp bound, a tractable strategy for constructing $M$-term approximations
is to chose them to be the $M$ largest values of $B(\pmb{\nu})$. Furthermore,
it is reasonable to assume that the bound 
$B(\pmb{\nu})$ is easier to compute than the coefficients themselves. 

For the Taylor and Legendre polynomials, sharp bounds on the coefficients are known.
Consider the Taylor series of a function $u$ given by 
$\sum_{\pmb{\nu} \in \mathcal{S}} c_{\pmb{\nu}} y^{\pmb{\nu}}$. Assuming that 
the function $u$ has certain smoothness assumptions outlined in
\cite[Proposition~1]{analysis_of_quasi_optimal} it is known that the Taylor coefficients
have the following bound,
	\begin{equation}
		\| c_{\pmb{\nu}} \| \le C \pmb{\rho}^{-\pmb{\nu}}.
	\end{equation}
In this case, the bounding function is $B(\pmb{\nu}) = C \pmb{\rho}^{-\pmb{\nu}}$.
On the other hand when $\{ \Psi_{\pmb{\nu}} \}$ is taken to be a tensor
product of Legendre Polynomials, \cite[Proposition~2]{analysis_of_quasi_optimal}
establishes the following bounds,
	\begin{equation}
		\| c_{\pmb{\nu}} \| \le C(\pmb{\rho}) \pmb{\rho}^{-\pmb{\nu}} 
		\Pi_{i=1}^d (2\nu_i - 1).
	\end{equation}
The estimates presented below are more general and apply to any polynomial
system for which there exists a bound $B(\bm{\nu})$ satisfying 
Assumption \ref{eq:bounds_assumption}.

\subsection{Useful Estimates}
We assume that $B(\pmb{\nu})$ takes the form 
	\begin{equation}
		B(\pmb{\nu}) = e^{-b(\pmb{\nu})}
	\end{equation}
where the function $b(\pmb{\nu})$ satisfies the following assumption.

\begin{assumption}{\cite[Assumption~3]{analysis_of_quasi_optimal}} \label{eq:bounds_assumption}
The map $b:[0,\infty)^d \rightarrow \mathbb{R}$ satisfies
	\begin{enumerate}
		\item $b(\pmb{0}) = 0$ and $b$ is continuous in $[0,\infty)^d$,
		
		\item $\frac{1}{\tau}b(\tau \pmb{\nu})$ is either increasing for $ \tau \in (0,\infty)$ and for all $\pmb{\nu} \in [0,\infty)^d$ 
		or decreasing for $ \tau \in (0,\infty)$ and for all $\pmb{\nu} \in [0,\infty)^d$,
		
		\item there exists $0 < c < C$ such that $c|\pmb{\nu}| < b(\pmb{\nu}) < C |\pmb{\nu}|$ as 
		$\pmb{\nu} \rightarrow \infty$ in the Euclidean norm.
	\end{enumerate}
\end{assumption}

According to \eqref{eq:upsilon_error}, if $u$ is a function whose expansion in the set of functions 
$\{ \Psi_{\pmb{\nu}} \}$ admits such a bound a convergence rate can be established by analyzing the series
	\begin{equation} \label{eq:series_convergence}
		\sum_{\pmb{\nu} \in \Lambda_M^c} e^{-b(\pmb{\nu})}.
	\end{equation}
The following result is proven in \cite{analysis_of_quasi_optimal} and exactly quantifies the 
convergence of the series \eqref{eq:series_convergence} when $\Lambda_M$ is taken to 
be the quasi-optimal index set.

\begin{lemma}{\cite[Theorem~2]{analysis_of_quasi_optimal}} \label{lem:main_quasi_approx}
	Consider the multi-indexed series 
	$\sum_{\pmb{\nu} \in \mathcal{S}} e^{-b(\pmb{\nu})}$ with $b:[0,\infty)^d \rightarrow \mathbb{R}$
	satisfying Assumption \ref{eq:bounds_assumption}. For $\tau \in (0,\infty)$, denote 
	$\mathcal{P}_\tau = \{ \pmb{\nu} \in [0,\infty)^d : b(\pmb{\nu}) \le \tau \}$ and 
	$\Lambda_M$ the set of indices corresponding to the $M$ largest 
	$e^{-b(\pmb{\nu})}$. Define $\mathcal{P}= \cap_{\tau \in \mathbb{R}^+}
	\left( \frac{1}{\tau} \mathcal{P}_\tau \right)$ 
	when $\frac{1}{\tau} b(\tau \pmb{\nu})$ is increasing 
	and
	$\mathcal{P} = \cup_{\tau \in \mathbb{R}^+} \left(\frac{1}{\tau} \mathcal{P}_\tau \right)$ 
	when $\frac{1}{\tau} b(\tau \pmb{\nu})$ is decreasing. If $\mathcal{P}$ is Jordan measurable,
	for any $\epsilon > 0$, there exists $M_\epsilon > 0$ depending on $\epsilon$
	such that 
	\begin{equation} \label{eq:quasi_approx_bounds}
		\sum_{\pmb{\nu} \not\in \Lambda_M} e^{-b(\pmb{\nu})} \le
		C_u(\epsilon) M 
		\exp \left ( - \left ( \frac{M}{|\mathcal{P}(1+\epsilon)|} \right )^{1/d} \right )	
	\end{equation}
	for all $M \ge M_\epsilon$. Here, $C_u(\epsilon) = (4e+4\epsilon e - 2)\frac{e}{e-1}$.
\end{lemma}

The error of our proposed network directly depends on \eqref{eq:quasi_approx_bounds}
and therefore we will require some estimates on $|\mathcal{P}|$. First, a useful lemma appearing in
\cite{analysis_of_quasi_optimal} shows that $|\mathcal{P}|$ is bounded and introduces a
useful characterization.

\begin{lemma}{\cite[Lemma~4]{analysis_of_quasi_optimal}} 
	Assume that $b:[0,\infty)^d \rightarrow \mathbb{R}$ satisfies Assumption \ref{eq:bounds_assumption}.
	Then, $0 < |\mathcal{P}| < \infty$. If $\mathcal{P}$ is Jordan measurable, there holds
		\begin{equation} \label{eq:P_limit}
			|\mathcal{P}| = \lim_{\tau \rightarrow \infty} \frac{1}{\tau^d} \cdot 
			\# \left ( \mathcal{P}_\tau \cap \mathbb{Z}^d \right ).
		\end{equation}
\end{lemma}

We establish a relationship between $|\mathcal{P}|$ and an integer $J$ in the following proposition
by considering the characterization of $|\mathcal{P}|$ given in \eqref{eq:P_limit} when the limiting 
index is assumed to be an integer.

\begin{proposition} \label{prop:J_bounds}
	For the set $\mathcal{P}$ as defined in Lemma \ref{lem:main_quasi_approx},  
	there exists $J \in \mathbb{N}$ (which may depend on $\epsilon$) such that
		\begin{equation}
			\left ( \frac{M}{|\mathcal{P}|(1+\epsilon)} \right )^{1/d} \le J \le \left ( \frac{2M}{|\mathcal{P}|} \right )^{1/d}.
		\end{equation}
\end{proposition}

The proof of this proposition is given in Appendix \ref{app:proof_J_bounds}.

\section{Construction and Analysis of the quasi-optimal Network} \label{sec:nn}
In this section, we will prove the main theorem showing the existence of a DNN 
which approximates arbitrarily well a quasi-optimal polynomial from a basis generated
by a tensor product of orthogonal polynomials. 

We will construct a Deep Neural Network $u_{NN}$ which uses ReLU activation functions. 
The \textit{complexity} of $u_{NN}$ 
is the number of nodes and edges in the graph induced by $u_{NN}$. 
The complexity is also the total number of weights and biases, 
since each edge is associated with a weight and each node is associated
with a bias. 
Denote by $u_{NN}^{(\Lambda_M)}$ the neural network which approximates 
the quasi-optimal approximation defined by the multi-index set $\Lambda_M$.
The network $u_{NN}^{(\Lambda_M)}$ is \textit{deep} in the sense it has many hidden layers
and its parameters and architecture will depend on the quasi-optimal index set. This dependence
also allows us to analyze the complexity of $u_{NN}^{(\Lambda_M)}$ defined as the total number of
weights and biases as well as its depth. In particular,
we will show that the complexity of the network is algebraic in $M$. 

\subsection{The quasi-optimal network} \label{sec:quasi_net}
In this section we prove our main result. The proof of this result constructs and analyzes a 
neural network which approximates a quasi-optimal polynomial approximation given both a 
polynomial basis generated from a tensor product of orthogonal polynomials and a quasi-optimal index set.
The construction of our network depends on being able to approximate the product of some inputs.
The work \cite{deep_learning_schwab} shows the existence of a neural network that approximates the 
product $\Pi_{i=1}^n x_i$.

\begin{lemma}{\cite[Corollary~3.3]{deep_learning_schwab}} \label{lem:product_nn}
	Let $\delta \in (0,1)$. 
	There exists a neural network denoted $\Tilde{\Pi}$ with $d$ input units such that for $x_1, \dots, x_d$ with
	$|x_i| \le 1$ for all $i$, it holds 
		\begin{equation}
			|\Pi_{j=1}^n x_j - \tilde{\Pi}(x_1,\dots,x_d)| \le \delta.
		\end{equation}
	Moreover, the complexity of the network $\tilde{\Pi}$, 
	i.e. the number of computational units and weights, is bounded by $C(1 + d\log(d/\delta))$ and
	the network $\tilde{\Pi}$ is no deeper than $C(1 + \log(d)\log(d/\delta))$ where $C = C(d,\delta)$ 
	is a constant depending on the number of
	inputs and the desired accuracy.
\end{lemma}

The construction of a network for approximating an orthogonal polynomial will use this 
product network. We now state 
our main theorem.

\begin{theorem} \label{thm:main_thm}
	Let $u: [0,1]^d \rightarrow \mathbb{R}$ satisfy Assumption \ref{as:function_of_interest}.
	Then, for any $M \in \mathbb{N}$ and $d \in \mathbb{N}$ there exists a neural network
	$u_{NN}(\bm{y})$ with $d$ inputs whose complexity 
	is bounded by $C M^{2/d +1}$, whose depth is bounded by 
	$1 + CM^{1/d} \log(M^{1/d})$ and which satisfies the error bound
		\begin{equation} \label{eq:network_error_subex}
			\left \| u(\pmb{y}) - u_{NN}(\pmb{y}) \right \|
			\le C M \exp \left( -\left( \frac{2M}{|\mathcal{P}|} \right)^{1/d} \right).
		\end{equation}
\end{theorem}

\begin{proof}
	Our goal is to estimate the error between a given function of interest $u$ and a network
	$u_{NN}$, which we will explicitly construct in this proof, and show that it achieves a 
	subexponential convergence rate in $M$. 
	Through our explicit construction of $u_{NN}$
	complexity bounds on the network are obtained. 
	We obtain the desired estimate by introducing an intermediate 
	approximation $u_Q$ which we choose to be a tensor product of orthogonal polynomials. 
	That is, for a given set of multi-dimensional orthogonal polynomials $\{ \Psi_{\pmb{\nu}} \}$, 
	the associated quasi-optimal polynomial is given by
	\begin{equation}
		u_{Q}(\pmb{y}) = \sum_{\pmb{\nu} \in \Lambda_M^{Qopt}} c_{\pmb{\nu}} \Psi_{\pmb{\nu}} (\pmb{y}),
	\end{equation}
	for a set of indices $\Lambda_M^{Qopt}$ as described in Section \ref{sec:quasi_optimal}.
	Notice that,
	\begin{equation} \label{eq:network_trick}
		\left \| u(\pmb{y}) - u_{NN}(\pmb{y}) \right \| \le
		\left \| u(\pmb{y}) - u_{Q}(\pmb{y}) \right \| + 
		\left \| u_{Q}(\pmb{y}) - u_{NN}(\pmb{y}) \right \|.
	\end{equation}
	Since we choose $u_Q$ to be a quasi-optimal polynomial approximation the first term
	can be bounded using Lemma \ref{lem:main_quasi_approx}. The second term is bounded by
	constructing a network $u_{NN}$ which approximates $u_Q$ with arbitrary accuracy. 
	The explicit construction of $u_{NN}$ will also reveal bounds on its complexity
	and depth. 
	 
	The network $u_{NN}$ will be composed of $M$ subnetworks $\tilde{\Psi}_{\pmb{\nu}}$
	which approximate the polynomials $\Psi_{\pmb{\nu}}$ for $\pmb{\nu} \in \Lambda_M^{Qopt}$.
	The network $\tilde{\Psi}_{\pmb{\nu}}$ will be constructed to compute the product 
	of $|\pmb{\nu}|_1$ numbers. The product of these numbers will approximate 
	the polynomial $\Psi_{\pmb{\nu}}$.
	Recall that any $d$-dimensional tensor product of a set of one-dimensional 
	orthogonal polynomials $\{ \psi_i \}_{i=1}^k$ may be written as 
		\begin{equation} \label{eq:tensor_product}
			\Psi_{\pmb{\nu}} (\pmb{y}) = \Pi_{i=1}^d \psi_{\nu_i}(y_i).
		\end{equation}
	The functions $\psi_{\nu_i}$ can be evaluated by a product of $\nu_i$ real numbers 
	by the fundamental theorem of algebra since they are assumed to have real roots, i.e.,
		\begin{equation} \label{eq:fundamental_theorem_express}
			\psi_{\nu_i} (y_i) = \Pi_{j=1}^{\nu_i} (y_i - r_j^{(\nu_i)}),
		\end{equation}
	where $\{r_j^{(\nu_i)}\}_{j=1}^{\nu_i}$ are the $\nu_i$ roots associated with the orthogonal polynomial 
	$\psi_{\nu_i}$ of degree $\nu_i$. Combining \eqref{eq:tensor_product} with 
	\eqref{eq:fundamental_theorem_express} we obtain
		\begin{equation} \label{eq:product_representation}
			\Psi_{\pmb{\nu}} (\pmb{y}) = \Pi_{i=1}^d \Pi_{j=1}^{\nu_i} (y_i - r_j^{(\nu_i)}).
		\end{equation}
	Notice that since $\psi_{\nu_i}$ are orthogonal on $[0,1]$, $r_j^{(\nu_i)} \in [0,1]$. Therefore,
	$|y_i - r_j^{(\nu_i)}|\le 1$ which is the assumption of Lemma \ref{lem:product_nn}.
	Therefore, $\Psi_{\pmb{\nu}}$ can be approximated by a network $\tilde{\Psi}_{\pmb{\nu}}$ 
	which first computes the values computes the product in \eqref{eq:product_representation}.
	
	The first layer of the network $\tilde{\Psi}_{\bm{\nu}}$ computes all the necessary  
	numbers of the form $(y_i - r_j^{(\nu_i)})$. 
	Notice that for any $y \in [0,1]$ and any $r \in [0,1]$,
		\begin{equation}
			(y-r) = \sigma(y-r) - \sigma(-(y-r))
		\end{equation}
	where $\sigma$ is the ReLU activation function. We can construct the appropriate inputs 
	to the subnetwork $\tilde{\Psi}_{\pmb{\nu}}$ using $2*|\pmb{\nu}|_1$ ReLU nodes. Alternatively,
	we can construct the first layer with $|\pmb{\nu}|_1$ nodes if we do not apply the ReLU activation
	function to each unit. After computing the necessary inputs the rest of the network $\tilde{\Psi}_{\bm{\nu}}$
	is composed a network which computes the product of these numbers. 
	By Lemma \ref{lem:product_nn},
	for every $\epsilon_{\pmb{\nu}} \in (0,1)$ there exists a network whose $|\pmb{\nu}|_1$
	inputs are $(y_i - r_j^{(\nu_i)})$ so that
		\begin{equation} \label{eq:network_poly_error}
			\| \Psi_{\pmb{\nu}} (\pmb{y}) - \tilde{\Psi}_{\pmb{\nu}} (\pmb{y}) \| < \epsilon_{\pmb{\nu}}.
		\end{equation}
	
	The function $u_Q$ can be approximated by a linear combination of the
	outputs of the networks $\tilde{\Psi}_{\pmb{\nu}}$, i.e.
		\begin{equation} \label{eq:network_output}
			u_Q(\pmb{y}) = \sum_{\nu \in \Lambda_M} c_{\pmb{\nu}} \Psi_{\pmb{\nu}} (\pmb{y}) \approx
			\sum_{\nu \in \Lambda_M} c_{\pmb{\nu}} \tilde{\Psi}_{\pmb{\nu}} (\pmb{y}) = u_{NN}(\pmb{y}).
		\end{equation}
	The desired network $u_{NN}$ can therefore be constructed so that it's output is the 
	final expression in \eqref{eq:network_output}.
	In light of \eqref{eq:network_output} we have 
		\begin{equation}
				\left \| u_{Q} - u_{NN}(\pmb{y}) \right \|
				\le 
				\sum_{\nu \in \Lambda_M} \| c_{\pmb{\nu}} \|
				\| \Psi_{\pmb{\nu}} (y) - \tilde{\Psi}_{\pmb{\nu}} \|.
		\end{equation}
	By applying the assumed coefficient bounds and \eqref{eq:network_poly_error} 
	we have,
		\begin{equation}
			 \sum_{\nu \in \Lambda_M} \| c_{\pmb{\nu}} \|
			\| \Psi_{\pmb{\nu}} (y) - \tilde{\Psi}_{\pmb{\nu}} \|
			\le \sum_{\nu \in \Lambda_M} e^{-b(\pmb{\nu})} \epsilon_{\pmb{\nu}}.
		\end{equation}
	The choice of approximation rate $\epsilon_{\bm{\nu}}$ of each of the networks 
	$\tilde{\Psi}_{\pmb{\nu}}$ has been left arbitrary up until this point
	but now we choose 
		\begin{equation} \label{eq:eps_nu}
			\epsilon_\nu = e^{b(\pmb{\nu}) - \left( \frac{2M}{|\mathcal{P}|} \right)^{1/d}}.
		\end{equation}
	Therefore, we have,
		\begin{equation} \label{eq:approx_network_error}
			\left \| u_{Q}(\pmb{y}) - u_{NN}(\pmb{y}) \right \|
			\le \sum_{\nu \in \Lambda_M} e^{-b(\pmb{\nu}) + b(\pmb{\nu}) 
			- \left( \frac{2M}{|\mathcal{P}|} \right)^{1/d}}
			= C M \exp \left( -\left( \frac{2M}{|\mathcal{P}|} \right)^{1/d} \right)
		\end{equation}
	as desired. 
	
	We will now justify our choice of $\epsilon_{\pmb{\nu}}$ as well as derive bounds on the depth and 
	complexity of the network $u_{NN}$. 
	In order for the result of Lemma \ref{lem:product_nn} to hold we require that $0 <\epsilon_{\pmb{\nu}} < 1$. Clearly,
	this is true for \eqref{eq:eps_nu} if and only if 
		\begin{equation} \label{eq:bounds_of_eps}
			- \infty < b(\pmb{\nu}) - \left( \frac{2M}{|\mathcal{P}|} \right)^{1/d} < 0.
		\end{equation}
	Recall that the index $\pmb{\nu}$ associated to $\epsilon_{\pmb{\nu}}$ is a member of the index set 
	$\Lambda_M$. Tthe characterization of this set given in \cite{analysis_of_quasi_optimal} is
		\begin{equation}
			\Lambda_M = \{ \nu \in [0,\infty): e^{-b(\nu)} \ge e^{-J} \} =  \{ \nu \in [0,\infty): b(\nu) \le J \}.
		\end{equation}
	Hence, we have
		\begin{equation} \label{eq:eps_bound_intermediate}
			e^{b(\pmb{\nu}) -  \left( \frac{2M}{|\mathcal{P}|} \right)^{1/d}} 
			\le e^{J - \left( \frac{2M}{|\mathcal{P}|} \right)^{1/d}}.
		\end{equation}
	Using Proposition \ref{prop:J_bounds} it is clear that 
		\begin{equation}
			e^{ -\left( \frac{2M}{|\mathcal{P}|} \right)^{1/d}} \le e^{-J}.
		\end{equation}
	Now combining this inequality with \eqref{eq:eps_bound_intermediate} we have
		\begin{equation}
			\epsilon_{\pmb{\nu}} = e^{b(\pmb{\nu}) -  \left( \frac{2M}{|\mathcal{P}|} \right)^{1/d}} 
			\le 1
		\end{equation}
	
	Having shown that we have chosen an admissible choice for $\epsilon_{\pmb{\nu}}$ we will 
	now be able to use the complexity estimates of Lemma \ref{lem:product_nn} to analyze the complexity
	of each network $\tilde{\Psi}_{\pmb{\nu}}$ in terms of $M$.
	According to \eqref{eq:product_representation}, the polynomial $\Psi_{\pmb{\nu}}$ 
	is a product of $|\nu|_1$ numbers,
	therefore the network $\tilde{\Psi}_{\pmb{\nu}}$ has $|\pmb{\nu}|_1$ inputs. 
	Recall that by Assumption \ref{eq:bounds_assumption} there exists $c > 0$ such that 
	$c |\pmb{\nu}|_1 \le b(\pmb{\nu})$. Then we can estimate the
	complexity and depth of $\tilde{\Psi}_{\pmb{\nu}}$ taking $\delta = \epsilon_{\pmb{\nu}}$. That is,
		\begin{equation}
			\begin{split}
			complexity(\tilde{\Psi}_{\pmb{\nu}}) & \le C(1 + |\pmb{\nu}|_1 \log(|\pmb{\nu}|_1)/\epsilon_{\pmb{\nu}}) \\
			& \le C( 1 + b(\pmb{\nu}) \log(b(\pmb{\nu})) - b(\pmb{\nu})\log(\epsilon_{\pmb{\nu}}) )\\
			& = C \left(1 + b(\pmb{\nu}) \left [ \log(b(\pmb{\nu})) - b(\pmb{\nu}) 
			+ \left( \frac{2M}{|\mathcal{P}|} \right)^{1/d} \right ] \right).
			\end{split}
		\end{equation}
	Noticing that $\log(b(\pmb{\nu})) - b(\pmb{\nu}) < 0$ since $b(\pmb{\nu}) > 0$ and
	recalling that $\pmb{\nu} \in \Lambda_M$ so that $b(\pmb{\nu}) \le J \le \left( \frac{2M}{|\mathcal{P}|} \right)^{1/d}$
	we have 
		\begin{equation}
			complexity(\tilde{\Psi}_{\pmb{\nu}}) \le C \left ( 1 + b(\pmb{\nu}) 
			\left( \frac{2M}{|\mathcal{P}|} \right)^{1/d} \right )
			\le C \left ( 1 + \left( \frac{2M}{|\mathcal{P}|} \right)^{2/d} \right ).
		\end{equation}
	A very similar calculation for the depth of $\tilde{\Psi}_{\pmb{\nu}}$ yields
		\begin{equation} \label{eq:depth_bound}
			depth(\tilde{\Psi}_{\pmb{\nu}}) \le C \left ( 1 + \left( \frac{2M}{|\mathcal{P}|} \right)^{1/d} 
			\log  \left( \frac{2M}{|\mathcal{P}|} \right)^{1/d}  \right ).
		\end{equation}
	
	Finally, we can estimate the complexity of the network $u_{NN}$ by summing the complexities of each 
	of the subnetworks that form it and analyze its depth by considering the deepest subnetwork $\tilde{\Psi}_{\pmb{\nu}}$. 
	First, we will consider the complexity of the first layer of the network. 
	Since each subnetwork $\tilde{\Psi}_{\pmb{\nu}}$ has $|\pmb{\nu}|_1$ inputs, 
	there are $\sum_{\nu \in \Lambda_M} |\pmb{\nu}|_1$ weights connecting the 
	inputs to the nodes on the first layer. Notice that 
		\begin{equation}
			\begin{split}
			\sum_{\pmb{\nu} \in \Lambda_M} |\pmb{\nu}|_1 
			& \le C \sum_{\pmb{\nu} \in \Lambda_M} b(\pmb{\nu}) \\
			& \le C \sum_{\pmb{\nu} \in \Lambda_M}  \left( \frac{2M}{|\mathcal{P}|} \right)^{1/d} \\
			& = C M \left( \frac{2M}{|\mathcal{P}|} \right)^{1/d} \\
			& = C M^{1/d + 1}
			\end{split}
		\end{equation}
	The outputs of each of the $\tilde{\Psi}_{\pmb{\nu}}$ blocks is connected to a single output
	node by $M$ connections whose weights are $c_{\pmb{\nu}}$. Therefore,
		\begin{equation}
			\begin{split}
				complexity(u_{NN}) & = C \left( M^{1/d+1} + \sum_{\pmb{\nu} \in \Lambda_M} 
				complexity(\tilde{\Psi}_{\pmb{\nu}}) + M \right ) \\
				& \le C \left ( M^{1/d+1} + M + M \left(1+ \left( \frac{2M}{|\mathcal{P}|} \right)^{2/d} \right )  \right )\\
				& \le C \left ( M^{1/d+1} + 2M +  \left( \frac{2}{|\mathcal{P}|} \right)^{2/d} M^{2/d+1} \right ) \\
				& \le C \left ( 1 + 2 + \left( \frac{2}{|\mathcal{P}|} \right)^{2/d}  \right ) M^{2/d+1} \\
				& \le C  ( 3 + 1) M^{2/d+1} \\
				& \le C M^{2/d +1}
			\end{split}
		\end{equation}
	since $ 2/|\mathcal{P}| < 1$ and $M \ge 1$.
	The depth of $u_{NN}$ is determined by deepest $\tilde{\Psi}_{\pmb{\nu}}$. Using the uniform depth bound given by 
	\eqref{eq:depth_bound}, one has
		\begin{equation}
			\begin{split}
			depth(u_{NN}) & = \max_{\pmb{\nu} \in \Lambda_M} depth(\tilde{\Psi}_{\pmb{\nu}}) \\
			& \le 1 + CM^{1/d} \log(C M^{1/d}) \\
			& \le 1 + CM^{1/d} \log(M^{1/d}).
			\end{split}
		\end{equation}
	Note that this estimate is not optimal since it relies on a uniform bound for all subnetworks 
	$\tilde{\Psi}_{\pmb{\nu}}$.
\end{proof}

\section{Conclusions}

We have shown that certain DNNs are capable of approximating functions with 
the same rate of approximation as quasi-optimal polynomial approximation. 
Our main result can be viewed as building a bridge
between traditional approximation problems, e.g. interpolation 
and polynomial approximation, and those considered in the machine learning
communities, e.g. classification and prediction. We hope that our work will direct 
attention to the connection between the classical approximations and 
machine learning. In particular, an approximation theoretic perspective may yield insight
into improved understanding of the generalization of neural networks. 
A direct application of the network constructed in this work would be to construct and initialize a 
network using a polynomial approximation of the training data. Such an approach 
may yield improved performance of neural networks used for approximating the solution 
to PDE.

\acks{We would like to Anton Dereventsov, Armenak Petrosyan, and Viktor Reshniak for many helpful discussions
during the formulation of this work.}

\bibliography{dnn_references}

\appendix

\section{Proof of Proposition \ref{prop:J_bounds}} \label{app:proof_J_bounds}

\begin{proof}
Let $ 1/2 > \epsilon > 0$, then according to \eqref{eq:P_limit} there exists an integer $N_\epsilon$ 
such that for all $J > N_\epsilon$,
	\begin{equation}
		\left | |\mathcal{P}| - \frac{1}{J^d} \cdot \# \left ( \mathcal{P}_J \cap \mathbb{Z}^d \right ) \right |
		\le \epsilon.
	\end{equation}
Then since $\epsilon < 1/2$ and $|\mathcal{P}| > 1$, we have 
	\begin{equation}
		- \epsilon |\mathcal{P}| < 
		|\mathcal{P}| - \frac{1}{J^d} \cdot \# \left ( \mathcal{P}_J \cap \mathbb{Z}^d \right ) 
		< \epsilon |\mathcal{P}| < \frac{1}{2} |\mathcal{P}|.
	\end{equation}
Therefore, for any $J > N_\epsilon$,
	\begin{equation}
			\left ( \frac{M}{|\mathcal{P}|(1+\epsilon)} \right )^{1/d} \le J \le \left ( \frac{2M}{|\mathcal{P}|} \right )^{1/d},
		\end{equation}
as desired.
\end{proof}

\end{document}